# GENERALIZED ITERATION, CATASTROPHES AND GENERALIZED SHARKOVSKY'S ORDERING

By Andrei Vieru

**Abstract**

We define iteration of functions that map $n$-dimensional vector spaces into $m$-dimensional vector spaces ($m$ at most equal to $n$). It happens that usual iteration and Fibonacci iterative methods become special cases of this 'generalized iteration'. Mathematical objects such as orbits, bifurcations, chaos, Feigenbaum constant, (generalized) Sharkovsky's ordering, (generalized) Julia and Mandelbrot sets and a new kind of catastrophe can be found and studied in this enlarged context.

**1. Generalized iteration of the first kind**

In spite of common ideas, iteration of maps don't suppose the domain and the range have the same number of dimensions.
A map $f: \mathbf{R}^n \to \mathbf{R}^m$ may very well be iterated, provided $m \leq n$ (otherwise we'll have an element of triviality in the definition).

Let's consider some arbitrary map $f: \mathbf{A}^n \to \mathbf{A}^m$ ($m \leq n$) and some arbitrary entry $(s_1, \ldots, s_n)$.

Writing $f(s_1, \ldots, s_n) = (g_1(s_1, \ldots, s_n), \ldots, g_m(s_1, \ldots, s_n)) = (s'_1, \ldots, s'_m)$, let us consider the sequence:

$u_1 = s_1$
$u_2 = s_2$
..........
$u_n = s_n$
$u_{n+1} = s'_1$
$u_{n+2} = s'_2$
.....
$u_{n+m} = s'_m$
$u_{n+m+1} = g_1(s_{m+1}, s_{m+2}, \ldots, s_{n-1}, s_n, s'_1, \ldots, s'_m) = s''_1$
$u_{n+m+2} = g_2(s_{m+2}, s_{m+3}, \ldots, s'_{m-1}, s'_m, s''_1) = s''_2$
.....
$u_{n+2m} = g_m(s_{2m}, s_{2m+1}, \ldots, s'_{m-1}, s'_m, \ldots, s''_{m-1}) = s''_m$
.....
$u_{n+km+j} = g_j(u_{km+j}, \ldots, u_{n+km+j-1})$
...........

This sequence, whose general term is $u_{n+km+j} = g_j(u_{km+j}, \ldots, u_{n+km+j-1})$ (for $k \geq 1$ and $j \geq 1$), may be viewed as a 'generalized F-iterative process'.

These rules constitute a generalization of the Fibonacci-type sequences. However they are **not** a generalization of the 'usual' iteration of maps $f: \mathbf{R}^n \to \mathbf{R}^n$

### 2. Generalized iteration of the second kind

Writing again $f(s_1,\ldots, s_n) = (g_1(s_1,\ldots, s_n),\ldots, g_m(s_1,\ldots, s_n)) = (s'_1,\ldots, s'_m)$, generalized iteration of the second kind (which is *also* a generalization of the Fibonacci-type iteration) is defined by the following rules:

$u_1 = s_1$
$u_2 = s_2$
……
$u_n = s_n$
$u_{n+1} = s'_1$
$u_{n+2} = s'_2$
……
$u_{n+m} = s'_m$
$u_{n+m+1} = g_1(s_{m+1}, s_{m+2},\ldots, s_{n-1}, s_n, s'_1,\ldots, s'_m) = s''_1$
$u_{n+m+2} = g_2(s_{m+1}, s_{m+2},\ldots, s_{n-1}, s_n, s'_1,\ldots, s'_m) = s''_2$
……
$u_{n+2m} = g_m(s_{m+1}, s_{m+2},\ldots, s_{n-1}, s_n, s'_1,\ldots, s'_m) = s''_m$
……
$u_{n+km+j} = g_j(u_{km+1},\ldots, u_{n+km})$
……

To make a distinction from the previously described *generalized F-iterative processes*, we'll call iteration based on these *rules generalized V-iterative process*. It's general term (for $k \geq 1$ and $j \geq 1$) is $u_{n+km+j} = g_j(u_{km+1},\ldots, u_{n+km})$

V-iteration *is* a generalization of 'usual' iteration. If $n = m$, we of course deal with normal iteration, because we may consider every subsequence $u_{kn+1}, u_{kn+2}\ldots, u_{kn+n}$ as one single point of $\mathbf{R}^n$.
In case $n = 2$, $m = 1$, $s_1 = 1$, $s_2 = 1$ and $f(x, y) = x + y$ we'll have the sequence of Fibonacci numbers, which is a special case of generalized F-iteration and V-iteration as well.

The concepts of *generalized F-iteration* and *generalized V-iteration* suggest us to formulate in higher dimension real or complex spaces all problems that already have been studied in lower dimension spaces, among which let us only mention those who lead to the concepts of Julia and Mandelbrot sets[1], to the discovery of the Feigenbaum

---

[1] One can construct, for example, generalized Julia and Mandelbrot sets based on the map $f: \mathbf{C}^2 \to \mathbf{C}$ with $f(z, w) = zw+c$ or $h(z, w) = z^2w^2+c$. Assuming $z=w$, generalized iteration of $f$

constant and of the Sharkovsky's ordering. For a given $\mathbf{C}^n \rightarrow \mathbf{C}^m$ mapping, we'll call *generalized Julia F-sets* (respectively *generalized Julia V-sets*) the subsets of $\mathbf{C}^n$ for which sequences generated by generalized F-iteration (respectively generalized V-iteration) are bounded. (Generalized Julia sets are not always themselves bounded!)

We formulate without proof
**A generalized F-iteration fixed point theorem**
If $y = f(x_1, x_2, \ldots, x_n)$ ($n \geq 2$) is a continuous map from $[0, 1]^n$ to $[0, 1]$ such that
$[\forall i \leq n \, [x_i = 0 \vee x_i = 1]] \Rightarrow f(x_1, x_2, \ldots, x_n) = 0$,
then it has fixed points, namely the roots of the equation $x = f(x, x, \ldots, x)$ that are situated in the interval $[0, 1]$.

### 3. Catastrophes of a new kind

Before coming close to some general conclusions, let's examine the family: $f_a(x, y) = ax(1-x)y(1-y)$ where $a$ is an arbitrary constant ($0 < a < 16$).

For sufficiently high $a$ values, we have 3 fixed points, namely the roots in $[0, 1]$ of the equation $a(x^4 - 2x^3 + x^2) - x = 0$, which is obtained from $ax(1-x)y(1-y) = x$, assuming $x = y$. Let $0 < r_1(a) < r_2(a)$ be these roots.
(For an arbitrary family of functions $ag(x, y)$: $[0, 1]^2 \rightarrow [0, 1]$, we'll have to consider the roots in $[0, 1]$ of the equation $ag(x, x) - x = 0$.)

For ***sufficiently high*** $a$ values, let's see what happens if we iterate $f_a(\alpha, \beta)$, constructing the sequence:
$\alpha, \quad \beta, \quad f_a(\alpha, \beta), \quad f_a(\beta, f_a(\alpha, \beta)), \quad f_a(f_a(\alpha, \beta), f_a(\beta, f_a(\alpha, \beta)))$, etc.

If the first and second arbitrarily chosen terms of the sequence – $\alpha$ and $\beta$ – correspond to a point $(\alpha, \beta)$ located in some open $\mathrm{A}_a$ area not too far from the border of the square $[0, 1]^2$, then the sequence converges to 0.

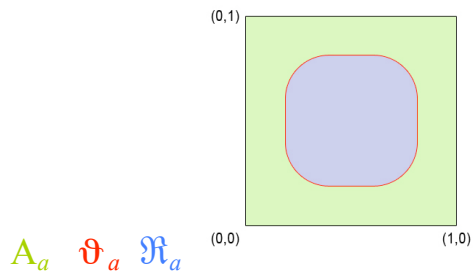

$\mathrm{A}_a \quad \vartheta_a \quad \mathfrak{R}_a$

---

engenders the sequence $z, z, z^2 + c, z^3 + zc + c, z^5 + 2z^3c + zc^2 + cz^2 + c^2 + c, \ldots$ whose features deserve to be studied, while 'classical iteration' of $g(z) = z^2 + c$ generates $z, z^2 + c, z^4 + 2cz^2 + c^2$, etc.

If ($\alpha$, $\beta$) is located on the $\vartheta_a$ boundary, then the sequence is attracted by $r_1(a)$ and $(r_1(a), r_1(a)) \in \vartheta_a$. We'll call $r_1(a)$ a **semi-instable** fixed point, in order to underline the fact that it attracts only the **one**-dimensional $\vartheta_a$ set of points in the square $[0, 1]^2$.

If ($\alpha$, $\beta$) is situated in the open $\Re_a$ area, then the sequence converges to $r_2(a)$.

A completely unexpected fact – for higher $a$ values – is not so much the appearance of an attracting orbit with period 3, but the dependence of the trifurcation point of the parameter $a$ on initial conditions, *i.e.* on initial values assigned to the variables $x$ and $y$. We'll call this kind of fixed points *semi-stable*, since they attract a **two**-dimensional subset of $[0, 1]^2$, but not the whole open area $\Re_a$.

For example, for $x = 0.6$ and $y = 0.3$, the trifurcation point is $a \approx 13.1986367…$, whereas for $x = 0.7$ and $y = 0.3$ the trifurcation point is $a \approx 12.782842211…$.

Moreover, trifurcation occurs as a sudden shift, **without any continuous growth** of the distance the elements of the orbit are separated by. Why for some range of $a$ values the attracting root $r_2(a)$ of the equation $a(x^2 + x^4 – 2x^3) – x = 0$ ceases to be a stable fixed point, without becoming 'completely instable' remains an open question: F-iterating $\mathbf{R}^n \to \mathbf{R}$ mappings, stable fixed points split into orbits with period $n+1$, but these orbits don't immediately attract *all* sequences (so we'll also call these orbits *semi-stable*): for some initial $x$ and $y$ values in $\Re_a$, sequences continue to converge to the (semi-stable) fixed point $r_2(a)$.

This new kind of *catastrophe* seems amazing: all ingredients of the process *are continuous*. Examining analogous $\mathbf{R}^n \to \mathbf{R}$ map families, we see that $n$ of the $n+1$ elements of the orbit obtained after the first '$(n+1)$-furcation' are always identical. The $n+1$-th element is located far from the others.

Example: the family $f_a(x, y) = ax(1–x)y(1–y)$ generates for $a = 12.782842211$, $x = 0.7$ and $y = 0.3$, a sequence that converges to $r_2(12.782842211) = 0.6541934769…$, while for $a = 12.782842212$, $x = 0.7$ and $y = 0.3$, the period 3 orbit is made of $0.7973063104767419…$, $0.515931143419749…$, $0.515931143419749…$

We see that if ($\eta_a$, $\theta_a$, $\theta_a$) is an attracting orbit, then $x = \eta_a$ and $y = \theta_a$ is a solution of the equation

$$a^4(x–x^2)^2(y–y^2)^3[1–a(x–x^2)(y–y^2)]^2\{1–a^2(y–y^2)^2(x–x^2)[1–a(x–x^2)(y–y^2)]\}–y = 0 \qquad (1)$$

obtained writing explicitly the beginning of the F-iterative sequence:
$u_1 = x$
$u_2 = y$
$u_3 = ax(1-x)y(1-y)$
$u_4 = a^2y^2(1-y)^2x(1-x)[1 - ax(1-x)y(1-y)]$
$u_5 = a^4x^2(1-x)^2y^3(1-y)^3[1-ax(1-x)y(1-y)]^2\{1-a^2y^2(1-y)^2x(1-x)[1-ax(1-x)y(1-y)]\}$

The equation **(1)** is the explicit form of the equation $u_2 = u_5$ (or $y=u_5$)

One can also see that $x = \theta_a$ and $y = \eta_a$ is a solution of the equation

$$a^2y(1-y)^2x(1-x)[1-ax(1-x)y(1-y)]-x=0 \qquad (2)$$

(i.e. of the equation $u_1 = u_4$ or $x = u_4$)

A lot of strange phenomena appear during calculations: for instance, choosing always $x = y$ and decreasing $x$ and $y$ values toward $r_1(a)$ we'll sometimes get a sequence converging to $r_2(a)$ and sometimes we'll get a sequence attracted by a period 3 orbit.

Before we reach the semi-instable fixed point $r_1(a)$, the sequence will change its attractor *infinitely many times*, as it does when we approach from the interior of $\mathfrak{R}_a$ any point of the curve $\vartheta_a$ following any direction.

Let us divide $\mathfrak{R}_a$ into two areas: the area $\mathfrak{I}_a \subset \mathfrak{R}_a$ made of the points attracted by the semi-stable fixed point $r_2(a)$, and the area $\wp_a \subset \mathfrak{R}_a$ made of the points attracted by the periodic orbit ($\eta_a$, $\theta_a$, $\theta_a$): $\mathfrak{I}_a$ and $\wp_a$ are separated by what we naturally call a *catastrophic boundary*. It is difficult to say where $\mathfrak{I}_a$ and $\wp_a$ are open and where they are not: coming closer to the boundary, whatever the direction we approach from, the convergence to the fixed point $r_2(a)$ and the 'convergence' to the period 3 orbit become as slow as one wishes! For $(x, y)$ initial values situated very close to the catastrophic boundary, we'll have to look at more and more terms to know whether the sequence will eventually converge to the fixed point or to the period 3 orbit.

12. 562691867… is nearly the lowest value of the parameter $a$ for which at least some initial values of $x$ and $y$ bring about trifurcation[2],

---

[2] namely $x = 0.5$ and $y = 0.5$

whereas 13.5 is the highest *a* value for which at least some *x* and *y* values still generate sequences converging to the semi-stable fixed point $r_2(13.5)=0.666666\ldots$ In fact, it is a trifurcation point: the last one.

As one can see, the function $z=f_{13.5}(x, y)=13.5x(1 - x)y(1 - y)$ has both *partial derivatives*[3] in point (0.6666…, 0.6666…) equal to –1. The three other points in the unit square that have the same trifurcation point – namely (1/3, 1/3) (2/3, 1/3) and (1/3, 2/3) have partial derivatives (1, 1), (–1, 1) and (1, –1) However, there are a lot of points (*x*, *y*) in which, for a parameter *a* value beyond its trifurcation point, partial derivatives are both > –1 or both < –1 or one smaller and the other greater than –1, etc.

For every given *a*, $\Im_a$ and $\wp_a$ are bi-dimensional sets; if we add the parameter *a* axis, we'll have to consider a 3-dimensional 'primitive' fractal structure. We believe it deserves to be studied.

For higher values of the parameter *a*, we haven't anymore trifurcations, but doubling-period process leading to chaos: after semi-stable and semi-instable fixed points we meet stable orbits with period 3, then 6, 12, 24 etc.

The critical values of the parameter *a* don't depend anymore on initial *x* and *y* values. $\Re_a$ becomes larger as the value of the parameter *a* increases. $A_a$ and $\vartheta_a$ are still attracted respectively by 0 and $r_1(a)$.

Computing the *a* values $\mu_1, \mu_2, \ldots, \mu_n, \ldots$ for which stable period $3\times 2^n$ orbits split into stable period $3\times 2^{n+1}$ orbits, then computing the $(\mu_n - \mu_{n-1})/(\mu_{n+1} - \mu_n)$ ratios, we quickly find numbers close to the Feigenbaum constant (4.669…), which is the limit of $(\mu_n - \mu_{n-1})/(\mu_{n+1} - \mu_n)$ when $n\to\infty$. Somewhere near $\mu_\infty \approx 13.78$ we reach chaos.

Remarkably, stable orbits appear again beyond the chaos point, as it happens in usual iteration of $[0, 1]\to[0, 1]$ maps.

Thus, within a large interval containing *a* = 13.97, we find a stable orbit with period 9 = 3×3 (and not at all of order 3, as it would have happened if we iterated an endomorphism of the interval), which soon splits into stable orbits with period 18, then 36, 72, 144, etc… leading again to chaos, just as it happens in usual iteration.

At about *a* = 13.883 we find a stable orbit with period 15 = 3×5 (and not anymore of order 5, as it would have happened if we iterated

---

[3] we propose the reader to introduce in the theory of generalized iteration some concept analogous to the concept of *eigenvalue* in classical theory of maps of the interval iteration. But, obviously, all geometric interpretations will be less comfortable to handle…

maps of the interval), which soon splits into stable orbits with period 30, then 60, 120, 240, etc… leading to chaos again, as it happens in usual iteration.

### 4. Generalized Sharkovsky's ordering

It is well known that in usual iteration of maps of the interval stable orbits appear before and beyond the chaos point in relationship with Sharkovsky's ordering. It is almost exactly what happens when we F-iterate continuous $\mathbf{R}^n \to \mathbf{R}$ or $[0, 1]^n \to [0, 1]$ maps! We only have to consider this ordering on $(n+1)\mathbf{N} \cup \{1\}$ instead of $\mathbf{N}$.

In our $[0, 1]^2 \to [0, 1]$ case, the Sharkovsky's ordering[4] in $3\mathbf{N} \cup \{1\}$ is the following one:

$3.3 < 3.5 < 3.7 < 3.9 < 3.11 < … < 3.2.3 < 3.2.5 < 3.2.7 < … < 3.2.2.3 < 3.2.2.5 < 3.2.2.7 < … < 3.2.2.2.3 < … < 3.2^4 < 3.2^3 < 3.2^2 < 3.2 < 3 < 1$

We conjecture that for some classes of smooth $\mathbf{R}^n \to \mathbf{R}$ or $[0, 1]^n \to [0, 1]$ function families[5], the ordering in which stable orbits periods appear before and beyond the chaos point is always related to the following one:

$3(n+1) < 5(n+1) < 7(n+1) < 9(n+1) < 11(n+1) < … < 2.3(n+1) < 2.5(n+1) < 2.7(n+1) < 2.9(n+1) < … < 2.2.3(n+1) < 2.2.5(n+1) < 2.2.7(n+1) <…< 2.2.2.3(n+1) < … < 2^4(n+1) < 2^3(n+1) < 2^2(n+1) < 2(n+1) < n+1 < 1$

We'll call it *generalized Sharkovsky's ordering* on $(n+1)\mathbf{N} \cup \{1\}$. It is yet not clear if this analogy permits to think about the existence of a *generalized Sharkovsky's theorem* in the context of *generalized F-iteration* of $\mathbf{R}^n \to \mathbf{R}$ or $[0, 1]^n \to [0, 1]$ maps[6].

---

[4] This ordering exists even if there might be no 'generalized Sharkovsky's theorem' for continuous $\mathbf{R}^n \to \mathbf{R}$ or $[0, 1]^n \to [0, 1]$ maps: such an ordering really occurs in the appearance of periodic orbits before and beyond chaos points when we F-iterate $\mathbf{R}^n \to \mathbf{R}$ ($n \geq 2$) maps.

[5] Namely, $[0, 1]^n \to [0, 1]$ map families for which if $\forall i \leq n$ $[x_i = 0 \vee x_i = 1]$ then $f(x_1, x_2, …, x_n) = 0$, that have one single maximum and for which some additional condition is satisfied (analogous to that of 'good Schwarzian derivative' required in the context of iteration of continuous maps interval). We believe that in the $\mathbf{R}^n \to \mathbf{R}$ or $[0, 1]^n \to [0, 1]$ cases, under reasonable conditions, stable periodic orbits are either fixed points or of order $m(n+1)$.

[6] It may turn out that the existence of periodic orbits of order $(n+1)m$ are *relevant* – in the sense that they imply the existence of orbits of all orders that follow in the 'generalized Sharkovsky's ordering', while periods of orders $(n+1)m-1$, $(n+1)m-$

Whether such a generalized Sharkovsky's theorem holds or not, period 3 implies, in the $\mathbf{R}^2 \rightarrow \mathbf{R}$ case, all but chaos!

In the $\mathbf{R}^3 \rightarrow \mathbf{R}$ case, when we increase parameter $a$ values, the map family $f_a(x, y, z) = ax(1-x)y(1-y)z(1-z)$ generates periodic orbits with periods of order 1, 4, 8, 16, 32, etc. leading to chaos[7].

Beyond chaos point, we find at $a = 54.51$ a stable orbit with period $20 = 4 \times 5$ that splits back to chaos through periods of order 40, 80, 160, etc.

Within a long interval around $a = 54.746$, we found a stable orbit of order $12 = 4 \times 3$. It splits and gets back to chaos through doubling-period process.

Stable orbits that appear before and beyond chaos points reached by generalized F- and V-iteration of $\mathbf{R}^n \rightarrow \mathbf{R}^m$ maps don't give birth to any Sharkovsky's ordering when $m \geq 2$. As we generally have $m$ parameters – so the chaos points are located on some curve inside the $\mathbf{R}^m$ space – the question doesn't really makes sense because we cannot consider 'increments' of the 'global parameters value' upon directions in the $\mathbf{R}^m$ space beyond the 'chaos curve': there is no suitable order there. In the $\mathbf{R}^3 \rightarrow \mathbf{R}^2$ case, the way from the very first fixed point (located in $\mathbf{R}^2$ and represented in the sequence by two consecutive terms) toward chaos is made – under some reasonable conditions – of usual doubling-period process.

All results mentioned above hold for other smooth functions with one single maximum inside $[0, 1] \times [0, 1]$ and such as $\forall x \in [0, 1]\ \forall y \in [0, 1]\ f(x, 0) = f(0, y) = 0$. (For instance, $\varphi_a(x, y) = a\sin(\pi x)\sin(\pi y)$.)

We also find that for function families like, for instance, $\phi_a(x, y) = a[x(1 - x) + y(1 - y)]$, which have one single maximum, but for which $\exists x\ \phi_a(x, 0) \neq 0 \wedge \exists y\ \phi_a(0, y) \neq 0$, the process that leads to chaos is also made of splitting orbits with periods 1, 3, 6, 12, 24, 48…

However, function families like $\phi_a(x, y)$ bring about catastrophes as trifurcation occurs, but they do not generate 'primitive fractal' sets: the $(x, x, x)$ diagonal crosses the $\varphi_a(x, y)$ functions graphs thrice (including the origin), while it crosses the $\phi_a(x, y)$ functions graphs only twice (including the origin). We believe that here resides the hidden reason for which primitive fractals appear or not.

Let us now consider a generalized F-iteration process based, as above, on $f_{a,b}(x, y, z) = (ag(x, y, z), bh(x, y, z))$. For a large class of

---

2,…, $(n+1)m-n$ are not – in the sense that they imply nothing. We propose the reader to shed light on this problem.

[7] All results concerning the appearance of the Feigenbaum constant, the catastrophic changes, the catastrophic boundaries and the 'primitive' fractal structure of the $\Im_a$ and $\wp_a$ sets can still be observed in this higher-dimensional case. As predicted, the initial fixed point splits directly into an order 4 periodic orbit.

mappings, fixing the value of one of the two parameters and increasing the other, we can compute the values for which bifurcations occur. Surprisingly, computing the limit ratio for consecutive distances between critical values of the 'free' parameter, we'll find again the Feigenbaum constant (4.669…).

Let $f(x_1, x_2,…, x_n) = (y_1, y_2,…, y_m)$ a function that maps $[0, 1]^n$ onto[8] $[0, 1]^m$. We'll consider that $f$ is made of $g_i(x_1, x_2,…, x_n) = y_i$ ($1 < i < m$).

We propose the reader to try to formulate and prove **a Generalized V-iteration fixed points theorem**[9].

We'll consider only some special cases:
If the system of equations
$g_i(x, x,…, x) = x$ had some solutions, these would have been fixed points of the form $(x, x,…, x)$. Generally, these solutions do not exist.
However, we may consider the system of equations with two unknowns:

$g_1(x_1, x_2, x_1, x_2,…, x_i, x_j) = x_i$      ($i = 1$ and $j = 2$ if $n$ is even; $i = 2$ and $j = 1$ if $n$ is odd)

$g_2(x_2, x_1, x_2, x_1, …, x_j, x_i) = x_j$

………

$g_m(x_k, x_l, x_k, x_l, …, x_p, x_r) = x_s$      ($s = k$ if $n$ is even and $s = l$ if $n$ is odd; $k = 1$ and $l = 2$ if $m$ is even; $k = 2$ and $l = 1$ if $m$ is odd; $p = k$ and $r = l$ if $n$ is even; $p = l$ and $r = k$ if $n$ is odd)

If this system admits solutions – which is a little bit more likely – then a pair of values that is a solution within the interval $[0, 1]$ will form a periodic orbit (in fact, a fixed point of the form $(x_1, x_2, x_1, x_2,…)$ in the mapped set $[0, 1]^m$). If this system has no solutions, we may still try to find solutions for the system of equations with three unknowns appearing, as previously, in a cycle[10]. And so on.

---

[8] we'll assume, for maps of the hypercube, that $\forall i\, [x_i=0 \lor x_i=1] \Rightarrow f(x_1, x_2,…, x_n)=0$

[9] In this context, we may consider as a fixed point any orbit of period $m$, because – we should not forget it – we consider the iteration of a $[0, 1]^n \to [0, 1]^m$ map. In order to avoid confusions, we even *must* do so, because it doesn't always make sense to consider, in $[0, 1]^m$ or $\mathbf{R}^m$, only fixed points of the form $(x, x,…, x)$, of the form $(x, y, x, y,…)$ or of the form $(x, y, z, x, y, z,…)$, etc.

[10] the system of equations will look like
$g_1(x_1, x_2, x_3, x_1, x_2, x_3,…, x_i, x_j, x_k)=x_i$ ($i=1\ j=2\ k=3$ if $n \equiv 0$ (mod. 3), $i=2\ j=3\ k=1$ if $n \equiv 3$, etc.)
$g_2(x_2, x_3, x_1, x_2, x_3,…, x_j, x_k, x_i)=x_j$ ($i=2\ j=2\ k=3$ if $n \equiv 0$ (mod. 3), $i=2\ j=3\ k=1$ if $n \equiv 3$, etc.)
……
$g_m(x_s, x_t, x_u, x_s, x_t, x_u,…, x_p, x_q, x_r)=x_w$ ($w=s$ if $n \equiv 0$ (mod. 3), $w=t$ if $n \equiv 1$ (mod. 3) $w=u$ if $n \equiv 2$ (mod. 3); $s=1\ t=2\ u=3$ if $m \equiv 1$ (mod. 3) $s=2\ t=3\ u=1$ if $m \equiv 2$ (mod. 3) $s=3\ t=1\ u=2$ if $m \equiv 0$ (mod.3); $p=s\ k=t\ r=u$ if $m \equiv 0$ (mod.3) $r=s\ p=t\ q=u$ if $m \equiv 1$ (mod. 3) $r=t\ p=u\ q=s$ if $m \equiv 2$ (mod.3)

A generalized V-iteration fixed points theorem may be formulated in an analogous way, taking care of writing the equations in a proper and suitable way.

V-iterating some $f_{a,b}: \mathbf{R}^3 \to \mathbf{R}^2$ map families, fixing one of the parameters and increasing the other, we quickly found the Feigenbaum constant. But perhaps we should not expect to find anymore the Feigenbaum constant considering the case when $a = b$, i.e. the case in which we have one single increasing parameter.

It is known that in the case of an area-preserving map $f_a: \mathbf{R}^2 \to \mathbf{R}^2$ the Feigenbaum constant value is 8,721097… Since 'generalized V-iteration' generalizes also 'normal' iteration – and not only the Fibonacci model – should we expect to find 8,721097… here again?

**Conclusion**

As the reader can see, a lot of interesting unsolved or not yet studied problems arise when studying generalized F- and V-iteration in the context of Dynamical Systems and Chaos Theory. However, the concepts of generalized F- and V-iteration are not specifically related to Dynamical Systems. (Usual iteration is not either.) It might be considered as combinatorial. As such, these concepts can and must be studied as well in completely different context. For instance in pure Algebra.

**Acknowledgments**
We express our gratitude to Dmitry Zotov for skillful computer programming and to Robert Vinograd for useful discussions on the subject.

**References**
(1) Vladimir Arnold *Catastrophe Theory*, 3rd ed. Berlin: Springer-Verlag, 1992.
(2) Bau-Sen Du *A collection of simple proofs of Sharkovky's theorem* arXiv:math/0703592v3
(3) Jonathan Borwein, David H. Bailey *Mathematics by Experiment: Plausible Reasoning in the 21st Century*
(4) Sharkovsky, A.N. [1964] : *Coexistence of cycles of continuous mapping of the line into itself*. Ukrainian Math. J., **16**, 61-71 [in Russian]; [1995] : Intern. J. Bifurcation and Chaos **5**, no. 5, 1263-1273.
(5) Andrei Vieru *Le gai Ecclésiaste,* Seuil 2007, p. 185-203

author's e-mail address: andreivieru@yahoo.fr